\documentclass[12pt]{amsart}
\usepackage{geometry}
\usepackage{color, enumerate,hyperref}
\usepackage{amssymb,amsmath,graphicx,pstricks}
\usepackage{amsmath,empheq}
\usepackage{bbm}

\makeatletter
\@namedef{subjclassname@2010}{%
  \textup{2010} Mathematics Subject Classification}
\makeatother

\newtheorem{thm}{Theorem}[section]

\theoremstyle{definition}

\newtheorem{theorem}{Theorem}

\newtheorem{rem}[thm]{Remark}

\newtheorem{defi}[thm]{Definition}

 \global\long\def\sbr#1{\left[ #1\right] }
 \global\long\def\cbr#1{\left\{  #1\right\}  }
 \global\long\def\rbr#1{\left(#1\right)}
 
 \global\long\def\E{\mathbb{E}}
 \global\long\def\P{\mathbb{P}}
 \global\long\def\R{\mathbb{R}}

 \global\long\def\TDD#1{}
 \global\long\def\dd#1{\textnormal{d}#1}

 \global\long\def\ab{[a,b]}
 \global\long\def\ra{\rightarrow}
 \global\long\def\TTV#1#2#3{\text{TV}^{#3}\!\rbr{#1,#2}}
 
 \global\long\def\UTV#1#2#3{\text{UTV}^{#3}\!\rbr{#1,#2}}
 \global\long\def\DTV#1#2#3{\text{DTV}^{#3}\!\rbr{#1,#2}}
 \global\long\def\ns{\infty}
 \global\long\def\f{:\left[a,b\right]\ra\R}

\newenvironment{dwd}{\par\noindent{\bf Proof.}}{\par\rightline{$\blacksquare$}}

\begin{document}

\baselineskip=17pt

\title{On a generalisation of the Banach Indicatrix Theorem}

\author[R. M. \L ochowski]{Rafa{\l } M. \L ochowski}
\address{{Department of Mathematics and Mathematical Economics}\\ {Warsaw School of Economics}\\{ul. Madali\'{n}skiego 6/8, 02-513 Warszawa, Poland}}
\address{{African Institute for Mathematical Sciences}\\{6-8 Melrose Road, Muizenberg 7945, South Africa}}
\email{rlocho314@gmail.com}

\date{}

\begin{abstract}
We prove that for any regulated function $f:\left[a,b\right]\rightarrow\mathbb{R}$ and $c\geq0,$
the infimum of the total variations of functions approximating $f$ with accuracy $c/2$ is equal$\int_{\mathbb{R}} n_{c}^{y} \dd y,$ where $n_{c}^{y}$ is the number of times that $f$ crosses the interval $[y,y+c].$ 
\end{abstract}



\maketitle

\section{Introduction, definitions, notation}

Let $a,$ $b$ such that $a<b$ be two fixed reals. The Banach Indicatrix
Theorem describes an interesting relationship between the total variation
of a continuous function $f:\left[a,b\right]\rightarrow\mathbb{R},$
defined as 
\[
\TTV f{\left[a,b\right]}{}:=\sup_{n}\sup_{a\leq t_{0}<t_{1}<\ldots<t_{n}\leq b}\sum_{i=1}^{n}\left|f\left(t_{i}\right)-f\left(t_{i-1}\right)\right|
\]
and the numbers of level crossings of the function $f,$ defined as
\[
N^{y}\left(f\right):=\#\left\{ x\in\left[a,b\right]:f\left(x\right)=y\right\} .
\]
($\#A$ denotes here the cardinality of the set $A$.) The function
$\mathbb{R}\ni y\mapsto N^{y}\left(f\right)\in\left\{ 0,1,2,\ldots\right\} \cup\left\{ +\infty\right\} $
is a Baire function of class $\leq 2$ and is called the \emph{Banach indicatrix}. The aforementioned relationship
has the form 
\[
\TTV f{\left[a,b\right]}{}=\int_{\mathbb{R}}N^{y}\left(f\right)\mathrm{d}y.
\]
It was proven by Banach in \cite{BanachIndicatrix:1925}, see \cite[Th\'{e}or\`{e}me 1.2]{BanachIndicatrix:1925}, and he is usually given full credit for it.  Vitali published the result (in the same journal) one year later, see \cite{Vitaliindicatrix:1926}, and Sergei M. Lozinskii \cite{Lozinskiindicatrix:1948}, \cite{Lozinskiindicatrix:1958}, generalized it to include the case of \emph{regulated} functions, i.e. functions $f:\left[a,b\right]\rightarrow\mathbb{R}$ admitting right limits $f\left(t+\right)$ at any point $t\in\left[a,b\right)$
and left limits $f\left(t-\right)$ at any point $t\in\left(a,b\right].$ For modern exposition of the proof of the Banach Indicatrix Theorem for continuous functions see \cite[Theorem 4.2.7]{BenedettoCzaja:2014} or \cite[Theorem (3.i)]{Cesari:1958}. 

Unfortunately, when $\TTV f{\left[a,b\right]}{}= +\ns$ this result seems to be useless. The purpose of this paper is to state and prove a meaningful generalisation of the Banach Indicatrix Theorem for \emph{any}
real regulated function, possibly with infinite total variation, in terms of \emph{segment crossings} rather than level crossings. 
Together with the numbers of segment and level crossings of the function $f,$
we will consider the numbers of segment and level upcrossings (or crossings from below) and numbers of segment and level downcrossings (or crossings from above). The  numbers of level upcrossings and downcrossings are related to the \emph{positive} and \emph{negative} variations of a real function:
\[
\UTV f{\left[a,b\right]}{}:=\sup_{n}\sup_{a\leq t_{0}<t_{1}<\ldots<t_{n}\leq b}\sum_{i=1}^{n}\left(f\left(t_{i}\right)-f\left(t_{i-1}\right)\right)_+,
\]
\[
\DTV f{\left[a,b\right]}{}:=\sup_{n}\sup_{a\leq t_{0}<t_{1}<\ldots<t_{n}\leq b}\sum_{i=1}^{n}\left(f\left(t_{i}\right)-f\left(t_{i-1}\right)\right)_-,
\]
where $x_+ = \max \left\{x,0 \right\},$ $x_- = \max \left\{-x,0 \right\}.$
Recall the Jordan decomposition:
\[
\TTV f{\left[a,b\right]}{} = \UTV f{\left[a,b\right]}{} + \DTV f{\left[a,b\right]}{}.
\]

It is possible to prove an analog of the Banach result for negative and positive variations and numbers of level upcrossings and downcrossings, see \cite[Section 9]{Cesari:1958}. 
To be more precise, let us start with several definitions. Let us define the number of
\emph{downcrossings} \emph{of $f$ from above the level $y+c$ to
the level $y.$} 
\begin{defi}\label{defd} Given a function $f:\left[a,b\right]\rightarrow\mathbb{R},$
for $c\geq0$ we put $\sigma_{0}^{c}=a$ and for $n=0,1,...$ 
\[
\tau_{n}^{c}=\inf\left\{ t>\sigma_{n}^{c}:t\leq b,f(t)>y+c\right\} ,\mbox{ }\sigma_{n+1}^{c}=\inf\left\{ t>\tau_{n}^{c}:t\leq b,f(t)<y\right\} .
\]
The \emph{number $d_{c}^{y}\left(f,\left[a,b\right]\right)$ of downcrossings
of $f$ from above the level $y+c$ to the level $y$ (or downcrossings by $f$ the closed segment $[y, y+c]$) on 
the interval }$\left[a,b\right]$ is defined as 
\[
d_{c}^{y}\left(f,\left[a,b\right]\right):=\max\left\{ n:\sigma_{n}^{c}\leq b\right\} .
\]
\end{defi} Similarly we define the number of \emph{upcrossings} \emph{of
$f$ from below the level $y$ to the level $y+c.$} 
\begin{defi}\label{defu}
Given a function $f:\left[a,b\right]\rightarrow\mathbb{R},$ for $c\geq0$
we put ${\tilde\sigma}_{0}^{c}=a$ and for $n=0,1,...$ 
\[
{\tilde\tau}_{n}^{c}=\inf\left\{ t>{\tilde\sigma}_{n}^{c}:t\leq b,f(t)<y\right\} ,\mbox{ }{\tilde\sigma}_{n+1}^{c}=\inf\left\{ t>{\tilde\tau}_{n}^{c}:t\leq b,f(t)>y+c\right\} .
\]
The \emph{number $u_{c}^{y}\left(f,\left[a,b\right]\right)$ of upcrossings
of $f$ from below the level $y$ to the level $y+c$ (or upcrossings by $f$ the closed segment $[y, y+c]$) on the interval $\left[a,b\right]$} is defined as 
\begin{equation}
u_{c}^{y}\left(f,\left[a,b\right]\right):=\max\left\{ n:{\tilde\sigma}_{n}^{c}\leq b\right\} .\label{eq:u_def}
\end{equation}
\end{defi} 
In all definitions we apply the convention that $\inf \emptyset = +\ns.$
At last, for $f$ and the interval $\left[a,b\right]$
as in two preceding definitions, we define the \emph{number $n_{c}^{y}\left(f,\left[a,b\right]\right)$
of crossings} \emph{by $f$ the segment $\left[y,y+c\right]$ (from
above or from below) on   the
interval }$\left[a,b\right]$ as 
\[
n_{c}^{y}\left(f,\left[a,b\right]\right):=d_{c}^{y}\left(f,\left[a,b\right]\right)+u_{c}^{y}\left(f,\left[a,b\right]\right).
\]
\begin{rem}
The precise definition of the number of level (up-, down-) crossings may be obtained by setting in the preceding definitions $c=0.$ We define the \emph{number $D^{y}\left(f,\left[a,b\right]\right)$ of crossings of $f$ the level $y$ from above on the interval} $\left[a,b\right]$ simply as 
\[
D^{y}\left(f,\left[a,b\right]\right):=d_{0}^{y}\left(f,\left[a,b\right]\right)
\]
and the \emph{number $U^{y}\left(f,\left[a,b\right]\right)$ of crossings
of $f$ the level $y$ from below on  
the interval }$\left[a,b\right]$ as 
\[
U^{y}\left(f,\left[a,b\right]\right):=u_{0}^{y}\left(f,\left[a,b\right]\right).
\]
Finally, we define the \emph{number $N^{y}\left(f,\left[a,b\right]\right)$
of crossings} \emph{of $f$ the level $y$ (from above or from below)
on   the interval} $\left[a,b\right]$
as 
\[
N^{y}\left(f,\left[a,b\right]\right):=D^{y}\left(f,\left[a,b\right]\right)+U^{y}\left(f,\left[a,b\right]\right).
\]
\end{rem}
\begin{rem} \label{indicatrix}
Let us notice that the number of level crossings  $N^{y}\left(f,\left[a,b\right]\right)$ introduced in the previous remark may differ from the Banach indicatrix $N^{y}\left(f\right)$ (when $f$ is continuous) or its generalisation for regulated $f,$ introduced by Lozinskii in \cite{Lozinskiindicatrix:1948}, which we will also denote by $N^{y}\left(f\right).$ However, it is not difficult to prove that the set $\left\{y \in \R: N^{y}\left(f\right) \neq N^{y}\left(f,\left[a,b\right]\right)\right\}$ is countable. Indeed, both numbers coincide for any real $y\notin f\rbr{\sbr{a, b}}$ and for any  $y \in f\rbr{\sbr{a, b}} \setminus \cbr{f(a), f(b)}$  which is not a local maximum or minimum of $f.$
It remains to prove that the set of local maxima and minima is countable. Indeed, for any local maximum $y$ there exist two rational numbers $o$ and $r>0$ such that $[o-r,o+r] \subset [a,b]$ and $y = \max_{t \in [o-r,o+r]} f(t).$ Since the mapping $y \mapsto \rbr{o,r}$ is injective the set of local maxima is countable. Similarly one proves that the set of local minima is countable. 

Our numbers $N^{y}\left(f,\left[a,b\right]\right),$ $U^{y}\left(f,\left[a,b\right]\right)$ and $D^{y}\left(f,\left[a,b\right]\right)$ coincide with Cesari's definitions of numbers of level crossings, upcrossings and downcrossings denoted by $N_e,$ $N_+$ and $N_-$ respectively, see \cite[p. 329]{Cesari:1958}.
\end{rem}

Naturally, the just defined numbers of segment or level crossings
may be infinite, however, if $f$ is regulated then for any $c>0$
and $y\in\R,$ $n_{c}^{y}\left(f,\left[a,b\right]\right)$ is a finite
number, cf. \cite[Theorem 2.1]{NorvaisaConcrete:2010}. The fact
that a regulated function on a compact interval crosses any non-degenerate segment only finitely
many times is closely related to the fact that for any regulated function
its truncated variation at the truncation level $c>0$ on 
a compact interval is always finite. For $f:\left[a,b\right]\ra\R,$
its {\em truncated variation} with the truncation parameter $c>0$ is defined as 
\begin{equation}
\TTV f{\left[a,b\right]}c:=\sup_{n}\sup_{a\leq t_{0}<t_{1}<\ldots<t_{n}\leq b}\sum_{i=1}^{n}\max\left\{ \left|f(t_{i})-f(t_{i-1})\right|-c,0\right\} .\label{eq:TVDefinition}
\end{equation}
Similarly, we define the {\em upward} and {\em downward truncated variations}
(with the truncation parameter $c>0$) respectively as 
\begin{equation}
\UTV f{\left[a,b\right]}c:=\sup_{n}\sup_{a\leq t_{0}<t_{1}<\ldots<t_{n}\leq b}\sum_{i=1}^{n}\max\left\{ f(t_{i})-f(t_{i-1})-c,0\right\} \label{eq:UTVDefinition}
\end{equation}
and 
\begin{equation}
\DTV f{\left[a,b\right]}c:=\sup_{n}\sup_{a\leq t_{0}<t_{1}<\ldots<t_{n}\leq b}\sum_{i=1}^{n}\max\left\{ f(t_{i-1})-f(t_{i})-c,0\right\} .\label{eq:DTVDefinition}
\end{equation}
It is possible to prove that $f\f$ is regulated iff $\TTV f{\left[a,b\right]}c <+\infty $ for any $c>0,$ see \cite[Fact 2.2]{LochowskiColloquium:2013}.
\begin{rem} \label{Rem_2}
For $g,h:\ab\ra\R$ let us denote 
\[\left\Vert g-h\right\Vert _{\ns}:=\sup_{t\in\ab}\left|g\left(t\right)-h\left(t\right)\right|.\]
The (downward-, upward-) truncated variation has an interesting 
variational property. It is possible to prove (see \cite[Theorem 4]{LochowskiGhomrasniMMAS:2015}) that $\TTV f{\left[a,b\right]}c$ is the \emph{attainable} infimum of total variations of functions uniformly approximating regulated $f$ with accuracy $c/2,$
\[
\TTV f{\left[a,b\right]}c = \inf\left\{\TTV g{\left[a,b\right]}{} \text{ where } g\f,\left\Vert f-g\right\Vert _{\ns} \leq c/2\right\}.
\]
Similarly, $\UTV f{\left[a,b\right]}c$ is the \emph{attainable} infimum of positive variations of functions uniformly approximating $f$ with accuracy $c/2,$
\[
\UTV f{\left[a,b\right]}c = \inf\left\{\UTV g{\left[a,b\right]}{} \text{ where } g\f,\left\Vert f-g\right\Vert _{\ns} \leq c/2\right\},
\]
and similarly
\[
\DTV f{\left[a,b\right]}c = \inf\left\{\DTV g{\left[a,b\right]}{} \text{ where } g\f,\left\Vert f-g\right\Vert _{\ns} \leq c/2\right\}.
\]
From \cite[Theorem 4]{LochowskiGhomrasniMMAS:2015} and the Jordan decomposition it also follows that 
\begin{equation} \label{Jordan}
\TTV f{\left[a,b\right]}c = \UTV f{\left[a,b\right]}c + \DTV f{\left[a,b\right]}c.
\end{equation}
\end{rem}
Now we are ready to state the main result of this article.

\begin{theorem} \label{main}
Let $f:\ab\ra\R$ be a regulated function. For any $c>0,$  the mappings $y \mapsto u_{c}^{y}\left(f,\left[a,b\right]\right)$  and $ y \mapsto d_{c}^{y}\left(f,\left[a,b\right]\right) $ are pointwise limits of step functions and 
the following equalities hold 
\begin{equation}
\UTV f{\left[a,b\right]}c=\int_{\R}u_{c}^{y}\left(f,\left[a,b\right]\right)\dd y,\label{eq:estim1}
\end{equation}
\begin{equation}
\DTV f{\left[a,b\right]}c=\int_{\R}d_{c}^{y}\left(f,\left[a,b\right]\right)\dd y\label{eq:estimdtv}
\end{equation}
and 
\begin{equation}
\TTV f{\left[a,b\right]}c=\int_{\R}n_{c}^{y}\left(f,\left[a,b\right]\right)\dd y.\label{eq:estimtv}
\end{equation}
\end{theorem} 
From this theorem, Remark \ref{indicatrix} and the classical monotone convergence theorem for the Lebesgue integral, letting $c\downarrow0$ we easily obtain Lozinskii's result as well as Cesari's result \cite[Section 9]{Cesari:1958} for positive and negative variations.

\begin{rem}
The natural question arises whether it is possible to obtain a further generalisation of Theorem \ref{main}, where some segments are distinguished from others and the following integral is considered
\[
\int_{\R} n_{c}^{y}\left(f,\left[a,b\right]\right)m(y)\dd y,
\]
where  $m:\R \ra [0, +\ns)$ is a non-negative, Borel-measurable density function. 
Similar generalisation for $N^y(f)$ is known in the form of the change of variables formula under minimal assumptions, see \cite{Hajlasz:1993}. 
\end{rem}

\section{Proof of Theorem 1}

The proof will go along similar lines as the proof of \cite[Theorem 8]{LochowskiGhomrasniMMAS:2015}. Another possible method of the proof of slightly weaker estimates for a c\`{a}dl\`{a}g $f$: 
\begin{equation*}
\UTV f{\left[a,b\right]}c=\int_{\R}u_{c}^{y}\left(f,\left[a,b\right]\right)\dd y \leq \UTV f{\left[a,b\right]}c + c,
\end{equation*}
\begin{equation*}
\DTV f{\left[a,b\right]}c=\int_{\R}d_{c}^{y}\left(f,\left[a,b\right]\right)\dd y \leq \DTV f{\left[a,b\right]}c + c
\end{equation*}
and 
\begin{equation*}
\TTV f{\left[a,b\right]}c=\int_{\R}n_{c}^{y}\left(f,\left[a,b\right]\right)\dd y \leq \TTV f{\left[a,b\right]}c + 2c
\end{equation*}
was outlined and then used  in \cite{LochowskiGhomrasni:2014} to prove limit theorems for numbers of segment crossings for diffusions and semimartingales.
(Recall that $f$ is c\`{a}dl\`{a}g  if it is right continuous and has left
limits.) This method utilized the fact that for any $c>0$ and any
starting value $x\in\left[f\left(a\right)-c/2,f\left(a\right)+c/2\right]$,
there exists a function $f^{c,x}\f,$ such that $\left\Vert f^{c,x}-f\right\Vert _{\ns}\leq c/2,$
$f^{c,x}\left(a\right)=x$ and $f^{c,x}$ has the minimal total variation
on any interval $\left[a,t\right],$ $t\in\left[a,b\right],$ among
all functions approximating $f$ with accuracy $c/2$ and attaining
the value $x$ at $a.$ This function coincides with the solution
of the so called Skorohod problem on $\left[-c/2,c/2\right]$ for
$f,$ and has an interesting property, namely that the number of times
that $f^{c,x}$ crosses (from above or from below) the level $y+c/2$
is almost the same as the number of times that $f$ crosses the segment
$\left[y,y+c\right],$ see \cite[Lemmas 3.3, 3.4]{LochowskiGhomrasni:2014}. For other problems similar to the Skorohod problem, like the play operator or taut strings see for example \cite{LochowskiGhomrasniMMAS:2015} or \cite{LifshitsSetterqvist:2015}.
\begin{dwd}
First, similarly as in \cite{LochowskiGhomrasniMMAS:2015}, we will prove Theorem \ref{main} for the family of step functions and then we will utilise the fact
that each regulated function is a uniform limit of step functions. 

\textbf{Step 1. Proof for step functions.} First we will assume that
$f$ has the representation 
\begin{equation*}
f\left(t\right)=\sum_{k=0}^{n}f_{2k}1_{\left\{ t\left(2k\right)\right\} }\left(t\right)+\sum_{k=0}^{n-1}f_{2k+1}1_{\left(t\left(2k+1\right);t\left(2k+2\right)\right)}\left(t\right),\label{eq:psi_representation}
\end{equation*}
where $a=t\left(0\right) =t(1) < t\left(2\right) = t(3) <\ldots< t(2n-2) = t(2n-1)<t\left(2n\right)=b.$

Let $f_{\rbr{0}} \leq  f_{\rbr{1}} \leq \ldots \leq f_{\rbr{2n}}$ be the non-decreasing rearrangement of the sequence $f_i,$ $i=0,1,\ldots, 2n.$ For any real $y_1, y_2$ such that $f_{\rbr{i_1-1}} < y_1 \leq y_2 < f_{\rbr{i_1}} $ and $f_{\rbr{i_2-1}} < y_1 + c \leq y_2 +c < f_{\rbr{i_2}} $ for some  $i_1, i_2 =1,\ldots, 2n,$ we have $u_{c}^{y_1}\left(f,\left[a,b\right]\right) = u_{c}^{y_2}\left(f,\left[a,b\right]\right)$ and $d_{c}^{y_1}\left(f,\left[a,b\right]\right) = d_{c}^{y_2}\left(f,\left[a,b\right]\right),$ thus $y \mapsto u_{c}^{y}\left(f,\left[a,b\right]\right)$ and $y \mapsto d_{c}^{y}\left(f,\left[a,b\right]\right)$ are step functions. 

Now we will prove (\ref{eq:estim1}).
The (upward-, downward-) truncated variations of $f$ and the numbers of segment (up-, down-) crossings of $f$ are equal to the (discrete versions of the) truncated variations and the numbers of segment crossings of the function $p:\left\{0,1,\ldots, 2n\right\} \ra \R$ defined as 
\[
p(i) = f_i \quad \text{for } i=0,1,\ldots, 2n.
\]
The truncated variation of $p$ is simply defined as
\[
\TTV p{\ab}c := \max_{m \leq 2n} \max_{0\leq i_0 < i_1 < \ldots < i_m \leq 2n} \sum_{j=1}^{m} \max \cbr{\left| p\rbr{i_j}  - p\rbr{i_{j-1}} \right|-c,0}
\]
and analogously are defined the upward- and downward- truncated variations of $p.$
The definitions of the numbers of segment up- and down- crossings of $p$ are obvious modifications of definitions \ref{defu} and \ref{defd}.

Thus, it is enough to prove the thesis for the (upward-, downward-) truncated variation of $p$ and the numbers of segment (up-, down-) crossings of $p.$

First, for $i=0,1,\ldots,2n-1$ we define 
\[
I_{U}\left(i\right)=\min\left\{ j\in\left\{ i+1,\ldots,2n\right\} :f_{j}>\min_{i\leq k<j}f_{k}+c\right\} ,
\]
\[
I_{D}\left(i\right)=\min\left\{ j\in\left\{ i+1,\ldots,2n\right\} :f_{j}<\max_{i\leq k<j}f_{k}-c\right\} 
\]
with the convention that $\min\emptyset=+\infty.$

Now we are ready to compare the upward truncated variation of $p$
with the integrated numbers of segment crossings. Assume that $I_{U}\left(0\right)\leq I_{D}\left(0\right)$ (the case $I_{U}\left(0\right)\geq I_{D}\left(0\right)$ is symmetric).

If $I_{U}\left(0\right)=+\infty$ then also $I_{D}\left(0\right)=+\infty$
and the function $p$ crosses no segment of the form $\left[y,y+c\right]$
and it has osscillation smaller or equal $c,$ 
\[
\left\Vert p\right\Vert _{osc}:=\max_{0\leq i<j\leq 2n}\left|p\left(j\right)-p\left(i\right)\right|\leq c,
\]
hence 
\[
\TTV p{[0,2n]}c=0
\]
and 
\[
n_{c}^{y}\left(p,\left[0,2n\right]\right)=0\quad\mbox{for all }y\in\R.
\]
Thus the thesis follows. 

Now assume that $I_{U}\left(0\right)<+\infty,$
thus $I_{U}\left(0\right)\leq2n,$ and define sequences $\left(I_{D,k}\right)_{k=-1}^{\infty},\mbox{ }\left(I_{U,k}\right)_{k=0}^{\infty}$
in the following way: $I_{D,-1}=0,$ $I_{U,0}=I_{U}\left(0\right)$
and for $k=0,1,...$ 
\[
I_{D,k}=\left\{ \begin{array}{ll}
I_{D}\left(I_{U,k}\right) & \text{if }I_{U,k}\leq2n-1,\\
+\infty & \mbox{otherwise},
\end{array}\right.I_{U,k+1}=\left\{ \begin{array}{ll}
I_{U}\left(I_{D,k}\right) & \text{if }I_{D,k}\leq2n-1,\\
+\infty & \mbox{otherwise}.
\end{array}\right.
\]
For $k=0,1,...,$ such that $I_{D,k-1}\leq2n$ let us denote \[m_{k}=\min_{I_{D,k-1}\leq j\leq\min\left\{ I_{U,k}-1,2n\right\} }f_{j}\]
and for $k=0,1,...,$ such that $I_{U,k}\leq2n$ let us denote \[M_{k}=\max_{I_{U,k}\leq j\leq\min\left\{ I_{D,k}-1,2n\right\} }f_{j}.\]

We see that on the interval $\left[I_{D,-1},\min\left\{ I_{D,0},2n\right\} \right]$
the function $p$ crosses the segment $\left[y,y+c\right]$ from below
exactly once iff $y\in\left(m_{0},M_{0}-c\right).$  It does not cross
the segment $\left[y,y+c\right]$ if $y\leq m_{0}$ or $y\geq M_{0}-c.$
See Figure 1. 
\begin{figure} \label{fig:m_0}
\centering
\includegraphics[width=9cm, angle = 270]{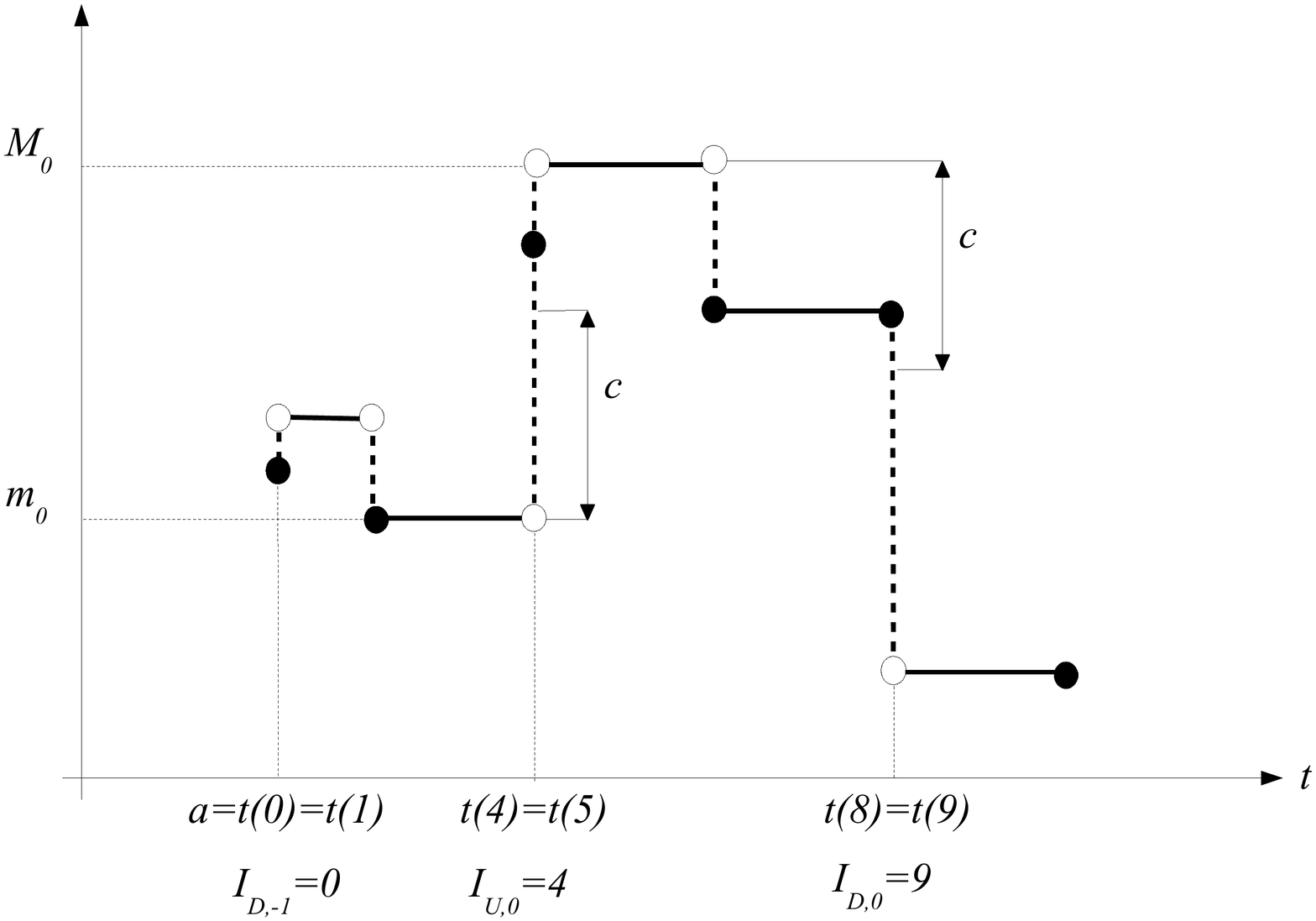}
\caption{Times $t\rbr{0}, t\rbr{1}, \ldots $ and values of $m_0$ and $M_0$ for a typical step path (the graph of the path is represented by the thick solid line and black dots).}
\end{figure}
Next, if $I_{D,0}\leq2n$ then the function $p$ crosses
no segment $\left[y,y+c\right]$ from below on the interval $\left[I_{D,0},\min\left\{ I_{U,1}-1,2n\right\}\right]$
but similarly as before, if $I_{U,1}\leq2n$ then it crosses the segment
$\left[y,y+c\right]$ from below exactly once on the interval $\left[I_{D,0},\min\left\{ I_{D,1},2n\right\} \right]$
iff $y\in\left(m_{1},M_{1}-c\right).$ It does not cross the segment
$\left[y,y+c\right]$ on $\left[I_{D,0},\min\left\{ I_{D,1},2n\right\} \right]$
if $y\leq m_{1}$ or $y\geq M_{1}-c.$ Continuing this reasoning in
the same way further we get that for any $k=0,1,\ldots$ if $I_{U,k}\leq2n$
then on the interval $\left[I_{D,k-1},\min\left\{ I_{D,k},2n\right\} \right]$
we have 
\begin{eqnarray*}
M_{k}-m_{k}-c & = & \int_{m_{k}}^{M_{k}-c}1\dd y=\int_{\R}u_{c}^{y}\left(p,\left[I_{D,k-1},\min\left\{ I_{D,k},2n\right\} \right]\right)\dd y
\end{eqnarray*}

Next, let us denote 
\[
K=\max\left\{ k\in\left\{0, 1,\ldots \right\}:I_{U,k}\leq 2n\right\} .
\]
We will prove that in fact 
\begin{equation}
\sum_{k=0}^{K}\left\{ M_{k}-m_{k}-c\right\} =\int_{\R}u_{c}^{y}\left(p,[0,2n]\right)\dd y.\label{eq:crossings_1}
\end{equation}
This may be easily proven using e.g. induction with respect to $K.$
From equality 
\[
M_{k}-m_{k}-c=\int_{\R}1_{\left(m_{k},M_{k}-c\right)}\left(y\right)\dd y
\]
it follows that to prove (\ref{eq:crossings_1}) it is sufficient
to prove that for any real $y$ 
\begin{equation}
u_{c}^{y}\left(p,\left[0,2n\right]\right)=\#\left\{ k\in\left\{ 0,1,\ldots,K\right\}:y\in\left(m_{k},M_{k}-c\right)\right\} .\label{eq:crossings_2}
\end{equation}
For $K=0$ this equality was already justified. Next, let $K>0$ and
for $y\in\R,$ let 
\[
k_{0}\left(y\right)=\max\left\{ k:\left(k=-1\right)\mbox{ or }\left(k\in\left\{ 0,1,\ldots,K-1\right\}\mbox{ and }y\in\left(m_{k},M_{k}-c\right)\right)\right\} .
\]
By the induction hypothesis (\ref{eq:crossings_2}) this means that there are no two indices
$I_{D,k_{0}\left(y\right)}\leq i<j\leq I_{D,K-1}-1$ such that 
\[
p\left(i\right)<y<p\left(j\right)-c.
\]
But this also means that 
\begin{eqnarray*}
u_{c}^{y}\left(p,\left[0,2n\right]\right) & = & u_{c}^{y}\left(p,\left[0,I_{D,K-1}\right]\right)+u_{c}^{y}\left(p,\left[I_{D,K-1},\min\left\{ I_{D,K},2n\right\}\right]\right)\\
 & = & \#\left\{ k\in\left\{ 0,1,\ldots,K-1\right\}:y\in\left(m_{k},M_{k}-c\right)\right\} +1_{\left(m_{k},M_{k}-c\right)}\left(y\right)\\
 & = & \#\left\{ k\in\left\{ 0,1,\ldots,K\right\}:y\in\left(m_{k},M_{k}-c\right)\right\} 
\end{eqnarray*}
which completes the inducion proof of (\ref{eq:crossings_1}).

Now we will prove the equality 
\begin{equation}
\sum_{k=0}^{K}\left\{ M_{k}-m_{k}-c\right\} =\UTV p{\left[0,2n\right]}c.\label{eq:utv_1}
\end{equation}
We define the following function $p^{c}:\left\{ 0,1,\ldots,2n\right\} \ra\R,$
\[
p^{c}\left(i\right)=\begin{cases}
m_{0}+\frac{c}{2} & \mbox{if }i<I_{U,0},\\
\max_{j\in\left\{ I_{U,k},I_{U,k}+1,\ldots,i\right\} }f_{j}-\frac{c}{2} & \mbox{if }I_{U,k}\leq i<I_{D,k}\mbox{ for some }k=0,1,\ldots,\\
\min_{j\in\left\{ I_{D,k},I_{D,k}+1,\ldots,i\right\} }f_{j}+\frac{c}{2} & \mbox{if }I_{D,k}\leq i<I_{U,k+1}\mbox{ for some }k=0,1,\ldots
\end{cases}.
\]
From the definition of $I_{U,k},$ $I_{D,k},$ $k=0,1,\ldots$ (and
the assumption $I_{U,0}<I_{D,0}$) it follows that the function $p^{c}$
approximates $p$ with accuracy $c/2,$ $\left\Vert p-p^{c}\right\Vert _{\ns}\leq c/2$
(see also \cite[ proof of Theorem 8]{LochowskiGhomrasniMMAS:2015}). From this and the equality
\[
p\left(j\right)-p\left(i\right)-\cbr{p\left(j\right)-p^{c}\left(j\right)}-\cbr{p^{c}\left(i\right)
-p\left(i\right)}=p^{c}\left(j\right)-p^{c}\left(i\right) \quad \text{for } i,j\in[0,2n], 
\]
it follows 
\[
\left(p\left(j\right)-p\left(i\right)-c\right)_{+}\leq\left(p^{c}\left(j\right)-p^{c}\left(i\right)\right)_{+} \quad \text{for } i,j\in[0,2n]
\]
and 
\begin{equation}
\UTV p{[0,2n]}c\leq\UTV{p^{c}}{[0,2n]}{}\label{eq:utv_2}
\end{equation}
(recall also Remark \ref{Rem_2}). Moreover, since for $k$ such that $I_{U,k}<+\ns,$
$p^{c}$ is non-decreasing on $\left[I_{U,k},\min\left\{ I_{D,k}-1,2n\right\} \right]$
and for $k$ such that $I_{D,k}<+\ns,$ $p^{c}$ is non-increasing
on $\left[I_{D,k},\min\left\{ I_{U,k+1}-1,2n\right\} \right]$
we see that 
\begin{equation}
\UTV{p^{c}}{[0,2n]}{}=\sum_{k=0}^{K}\left\{ M_{k}-m_{k}-c\right\} .\label{eq:utv_3}
\end{equation}
On the other hand, since $\UTV p{\cdot}c$ is a super-additive function
of the interval, i.e. for $0<i<2n$ 
\[
\UTV p{[0,2n]}c\geq\UTV p{\left[0,i\right]}c+\UTV p{\left[i,2n\right]}c,
\]
we have that 
\begin{eqnarray}
\UTV p{[0,2n]}c & \geq & \sum_{k=0}^{K-1}\UTV p{\left[I_{D,k-1},I_{D,k}\right]}c+\UTV p{\left[I_{D,K-1},2n\right]}c\nonumber \\
 & \geq & \sum_{k=0}^{K}\left\{ M_{k}-m_{k}-c\right\} .\label{eq:utv_4}
\end{eqnarray}
Now, from (\ref{eq:utv_2})-(\ref{eq:utv_4}) we get (\ref{eq:utv_1}).
Finally, from (\ref{eq:crossings_1}) and (\ref{eq:utv_1}) we obtain
the equality 
\[
\UTV p{[0,2n]}c=\int_{\R}u_{c}^{y}\left(p,\left[0,2n\right]\right)\dd y
\]
which may be translated into the equality
\[
\UTV f{\ab}c=\int_{\R}u_{c}^{y}\left(f,\left[a,b\right]\right)\dd y.
\]

Similarly, considering segment downcrossings and the downward truncated
variation of $p$ on the intervals $\left[I_{U,k},I_{U,k+1}\right],$
$k=0,1,\ldots,$ one obtains 
\[
\DTV p{[0,2n]}c=\int_{\R}d_{c}^{y}\left(p,\left[0,2n\right]\right)\dd y
\]
which translates into the equality
\[
\DTV f{\ab}c=\int_{\R}d_{c}^{y}\left(f,\left[a,b\right]\right)\dd y.
\]

Finally, from the definition of $n_{c}^{y}$ and the equality $\TTV f{\ab}c=\UTV f{\ab}c+\DTV f{\ab}c$ (recall (\ref{Jordan}) in Remark \ref{Rem_2})
we get 
\[
\TTV f{\ab}c=\int_{\R}n_{c}^{y}\left(f,\left[a,b\right]\right)\dd y.
\]

\textbf{Step 2. Proof for arbitrary regulated functions.} 
We will use Lemma 21 from \cite{LochowskiGhomrasniMMAS:2015} and the fact that any regulated function $f\f$ is a uniform limit of step functions, see e.g. \cite[Theorem 7.6.1]{Dieudonne:1969}.

Let $\varepsilon\in\left(0,c/2\right)$ and
$f^{\varepsilon}\f$ be a step function such that $\left\Vert f-f^{\varepsilon}\right\Vert _{\ns}\leq\varepsilon.$
From the very definition of segment crossings from below we get that
\begin{equation}
u_{c}^{y}\left(f,\ab\right)\leq u_{c-2\varepsilon}^{y+\varepsilon}\left(f^{\varepsilon},\ab\right).\label{eq:level_1}
\end{equation}
Indeed, each upcrossings of $f$ from below the level $y$ to the
level $y+c$ corresponds to two times $s,t$ such that $a\leq s<t\leq b$
and $f\left(s\right)<y$ and $f\left(t\right)>y+c.$ From this and
$\left\Vert f-f^{\varepsilon}\right\Vert _{\ns}\leq\varepsilon$ we
immediately get 
\[
f^{\varepsilon}\left(s\right)\leq f\left(s\right)+\varepsilon<y+\varepsilon\mbox{ and }f^{\varepsilon}\left(t\right)\geq f\left(t\right)-\varepsilon>y+c-\varepsilon=\left(y+\varepsilon\right)+c-2\varepsilon,
\]
thus times $s,t$ also correspond to upcrossings of $f^{\varepsilon}$
from below the level $y+\varepsilon$ to the level $\left(y+\varepsilon\right)+c-2\varepsilon.$
Similarly, considering upcrossings of $f^{\varepsilon}$ from below
the level $y-\varepsilon$ to the level $\left(y-\varepsilon\right)+c+2\varepsilon$
one proves that 
\begin{equation}
u_{c+2\varepsilon}^{y-\varepsilon}\left(f^{\varepsilon},\ab\right)\leq u_{c}^{y}\left(f,\ab\right).\label{eq:level_2}
\end{equation}
Sending $\varepsilon$ to $0$ we see that $y \mapsto u_{c}^{y}\left(f,\ab\right)$ is a pointwise limit of step functions (since $y \mapsto u_{c-2\varepsilon}^{y+\varepsilon}\left(f^{\varepsilon},\ab\right)$ and $y \mapsto u_{c+2\varepsilon}^{y-\varepsilon}\left(f^{\varepsilon},\ab\right)$ are step functions). 
Now, from (\ref{eq:level_1}), (\ref{eq:level_2}) and Step 1 for step
functions we get 
\[
\UTV{f^{\varepsilon}}{\ab}{c+2\varepsilon}\leq\int_{\R}u_{c}^{y}\left(f,\ab\right)\dd y\leq\UTV{f^{\varepsilon}}{\ab}{c-2\varepsilon}.
\]
Now, sending $\varepsilon$ to $0$ and applying \cite[Lemma
21]{LochowskiGhomrasniMMAS:2015} we get the equality 
\[
\UTV f{\ab}c=\int_{\R}u_{c}^{y}\left(f,\ab\right)\dd y.
\]
Analogous equalities (\ref{eq:estimdtv}) and (\ref{eq:estimtv}) for $\DTV f{\ab}c$ and $\TTV f{\ab}c$ may be
justified in the same way. 

\end{dwd}

\section{Example of application}

In this section we will give an example of application of Theorem \ref{main}. Let $B_{t},$ $t\geq0,$ be a standard Brownian motion and $W_{t}=B_{t}+\mu t,$
$t\geq0$ be a standard Brownian motion with drift $\mu.$ In \cite{LochowskiSPA:2011} the following formula for the Laplace transform of the function $t\mapsto\E\UTV W{\left[0,t\right]}c$
was given \cite[formula (25)]{LochowskiSPA:2011}: if $\nu$ is a complex
number with negative real part, $\Re\left(\nu\right)<0,$ then 
\begin{equation}
\int_{0}^{+\ns}e^{\nu t}\E\UTV W{\left[0,t\right]}c\dd t=\frac{e^{\mu c}\sqrt{\mu^{2}-2\nu}}{2\nu^{2}\sinh\left(c\sqrt{\mu^{2}-2\nu}\right)}.\label{eq:Lapl_tr}
\end{equation}
Using Theorem \ref{main} we will prove this formula when $v=-\nu$ is a positive
real number. Notice that if $v=-\nu$ is a positive real number and
$\tau$ is an exponential random variable with mean $1/v,$ independent
from $B,$ then 
\[
\int_{0}^{+\ns}e^{\nu t}\E\UTV W{\left[0,t\right]}c\dd t=\frac{1}{v}\E\UTV W{\left[0,\tau\right]}c.
\]
Thus, to prove formula (\ref{eq:Lapl_tr}) it is sufficient to prove
that 
\[
\E\UTV W{\left[0,\tau\right]}c=\frac{e^{\mu c}\sqrt{\mu^{2}+2v}}{2v\sinh\left(c\sqrt{\mu^{2}+2v}\right)}.
\]
First, for $y\in\R$ let us calculate $\E u_{c}^{y}\left(W,\left[0,\tau\right]\right).$
For $y\in\R$ let $\tau_{y}$ be the first hitting time of $y$
by $W,$ i.e. 
\[
\tau_{y}=\inf\left\{ t\geq0:W_{t}=y\right\} .
\]
Let $y\geq0,$ the event $u_{c}^{y}\left(W,\left[0,\tau\right]\right)\geq1$
is equivalent with the event $\left\{ \tau_{y+c}\leq\tau\right\} $
and by \cite[formula 1.1.2, p. 250]{BorodinSalminen:2002} we have 
\begin{eqnarray}
\P\left(u_{c}^{y}\left(W,\left[0,\tau\right]\right)\geq1\right) & = & \P\left(\tau_{y+c}\leq\tau\right)\nonumber \\
 & = & \P\left(\sup_{0\leq s\leq\tau}W_{s}\geq y+c\right)=e^{\mu\left(y+c\right)-\left(y+c\right)\sqrt{\mu^{2}+2v}}.\label{eq:Lapl1}
\end{eqnarray}
Now, let $y<0.$ By \cite[formula 1.2.2, p. 251]{BorodinSalminen:2002} we
have
\[
\P\left(\tau_{y}\leq\tau\right)=\P\left(\inf_{0\leq s\leq\tau}W_{s}\leq y\right)=e^{\mu y+y\sqrt{\mu^{2}+2v}}.
\]
Next, for $y<0$ the event $u_{c}^{y}\left(W,\left[0,\tau\right]\right)\geq1$
is equivalent with the event $\left\{ \tau_{y}\leq\tau\right\} \cap\left\{ \sup_{\tau_{y}\leq s\leq\tau}\left(W_{s}-W_{\tau_{y}}\right)\geq c\right\} $
and by the strong Markov property of the Brownian motion and the lack
of memory of the exponential distribution we get
\[
\P\left(\sup_{\tau_{y}\leq s\leq\tau}\left(W_{s}-W_{\tau_{y}}\right)\geq c|\tau_{y}\leq\tau\right)=\P\left(\sup_{0\leq s\leq\tau}W_{s}\geq c\right)
\]
thus
\begin{eqnarray}
\P\left(u_{c}^{y}\left(W,\left[0,\tau\right]\right)\geq1\right) & = & \P\left(\tau_{y}\leq\tau\right)\P\left(\tau_{c}\leq\tau\right)\nonumber \\
 & = & e^{\mu y+y\sqrt{\mu^{2}+2v}}e^{\mu c-c\sqrt{\mu^{2}+2v}}\nonumber \\
 & = & e^{\mu\left(y+c\right)+\left(y-c\right)\sqrt{\mu^{2}+2v}}.\label{eq:Lapl2}
\end{eqnarray}
Now, for $y\in\R$ let us define the following stopping time 
\[
\upsilon_{y}=\inf\left\{ t\geq0:u_{c}^{y}\left(W,\left[0,t\right]\right)\geq1\right\} .
\]
By the strong Markov property of the Brownian motion and the lack
of memory of the exponential distribution we get 
\begin{eqnarray*}
\P\left(u_{c}^{y}\left(W,\left[0,\tau\right]\right)\geq2\right) & = & \P\left(\upsilon_{y}\leq\tau\right)\P\left(u_{c}^{y}\left(W,\left[\upsilon_{y},\tau\right]\right)\geq1|\upsilon_{y}\leq\tau\right)\\
 & = & \P\left(u_{c}^{y}\left(W,\left[0,\tau\right]\right)\geq1\right)\P\left(u_{c}^{-c}\left(W,\left[0,\tau\right]\right)\geq1\right)\\
 & = & \P\left(u_{c}^{y}\left(W,\left[0,\tau\right]\right)\geq1\right)e^{-2c\sqrt{\mu^{2}+2v}}.
\end{eqnarray*}
Similarly, for all $n=1,2,3,\ldots,$ 
\begin{equation}
\P\left(u_{c}^{y}\left(W,\left[0,\tau\right]\right)\geq n\right)=\P\left(u_{c}^{y}\left(W,\left[0,\tau\right]\right)\geq1\right)e^{-2c\left(n-1\right)\sqrt{\mu^{2}+2v}}.\label{eq:Lapl3}
\end{equation}
From (\ref{eq:Lapl3}) we get 
\begin{eqnarray}
\E u_{c}^{y}\left(W,\left[0,\tau\right]\right) & = & \P\left(u_{c}^{y}\left(W,\left[0,\tau\right]\right)\geq1\right)\sum_{n=1}^{\ns}e^{-2c\left(n-1\right)\sqrt{\mu^{2}+2v}}\nonumber \\
 & = & \P\left(u_{c}^{y}\left(W,\left[0,\tau\right]\right)\geq1\right)\frac{1}{1-e^{-2c\sqrt{\mu^{2}+2v}}}.\label{eq:Lapl4}
\end{eqnarray}
Finally, from Theorem 1, (\ref{eq:Lapl1}), (\ref{eq:Lapl2}) and
(\ref{eq:Lapl4}) we have 
\begin{eqnarray*}
\E\UTV W{\left[0,\tau\right]}c & = & \E\int_{\R}u_{c}^{y}\left(W,\left[0,\tau\right]\right)\dd y=\int_{\R}\E u_{c}^{y}\left(W,\left[0,\tau\right]\right)\dd y\\
 & = & \int_{0}^{+\ns}\frac{e^{\mu\left(y+c\right)-\left(y+c\right)\sqrt{\mu^{2}+2v}}}{1-e^{-2c\sqrt{\mu^{2}+2v}}}\dd y \\
&& \quad +\int_{-\ns}^{0}\frac{e^{\mu\left(y+c\right)+\left(y-c\right)\sqrt{\mu^{2}+2v}}}{1-e^{-2c\sqrt{\mu^{2}+2v}}}\dd y\\
 & = & \frac{e^{\mu c-c\sqrt{\mu^{2}+2v}}}{1-e^{-2c\sqrt{\mu^{2}+2v}}}\left(\frac{1}{\sqrt{\mu^{2}+2v}-\mu}+\frac{1}{\sqrt{\mu^{2}+2v}+\mu}\right)\\
 & = & \frac{e^{\mu c}}{e^{c\sqrt{\mu^{2}+2v}}-e^{-c\sqrt{\mu^{2}+2v}}}\frac{2\sqrt{\mu^{2}+2v}}{2v} \\
& = &\frac{e^{\mu c}\sqrt{\mu^{2}+2v}}{2v\sinh\left(c\sqrt{\mu^{2}+2v}\right)}.
\end{eqnarray*}

\end{document}